\documentclass[11pt]{amsart}
\newtheorem{theo}{Theorem}[section]
\newtheorem{lemma}[theo]{Lemma}

\theoremstyle{definition}
\newtheorem{defi}[theo]{Definition}

\begin{document}

\title{Families of K3 surfaces}
\author[R.E.Borcherds]{Richard E. Borcherds$^\dagger$}
\thanks{$(\dagger)$ Partially supported by Grant DMS-9401186 and by a 
Royal Society professorship}
\address{Mathematics department, Evans Hall, 3840, UC Berkeley, CA 94720-3840}
\curraddr{D.P.M.M.S., 16 Mill Lane, Cambridge CB2 1SB, UK}
\email{reb@pmms.cam.ac.uk}
  
\author[L.Katzarkov]{Ludmil Katzarkov}
\address{M.S.R.I., 1000 Centennial Drive, Berkeley, CA 94720} 
\curraddr{Department of Mathematics, UC Irvine, Irvine, CA 92697-3875}
\email{lkatzark@math.uci.edu}

\author[T.Pantev]{Tony Pantev$^\ddagger$}
\address{Department of Mathematics, MIT, Cambridge, MA 02139}
\email{pantev@math.mit.edu}
\thanks{$(\ddagger)$ Partially supported by NSF grant DMS-9500712.}

\author[N.I.Shepherd-Barron]{N.~I.~Shepherd-Barron$^{*}$}
\address{D.P.M.M.S., 16 Mill Lane, Cambridge CB2 1SB, UK}
\email{nisb@pmms.cam.ac.uk}
\thanks{$(*)$ Partially supported by a European Community HCM grant.}



\maketitle

\tableofcontents

\section {Introduction}

We will prove the following theorem and give some examples to show
that most of the conditions in it are necessary.  Recall that a family
$X\to B$ of varieties over a base space $B$ is called isotrivial if
there is an \'etale covering $\tilde B\to B$ such that
$X\times_B\tilde B$ is a trivial family over $\tilde B$.

\begin{theo}
Any complete family of minimal K\"ahler surfaces 
of Kodaira dimension 0 and constant Picard number is isotrivial.
\end{theo}

This is a generalization to surfaces of the well known fact that any
complete family of complex elliptic curves is isotrivial, because any
complex elliptic curve is automatically minimal, K\"ahler, of Kodaira
dimension 0, and has Picard number 1. Roughly speaking, it gives some
cases when moduli spaces of surfaces contain no complete subvarieties.

In this paper  all varieties are smooth, and are defined over
$\mathbb{C}$ except when we explicitly state otherwise. 

Recall that any minimal K\"ahler surface of Kodaira dimension 0 is
abelian, hyperelliptic, Enriques, or K3.  We prove theorem 1.1 by
treating these cases separately.  We can quickly
dispose of most of these cases and reduce to the case of
projective K3 surfaces
as follows.  

We often want to deduce that a compact family of surfaces with extra
structure (including at least a polarization) is
isotrivial from the fact that a suitable moduli space is quasiaffine
(so that all fibers in the family are isomorphic as the image of any
compact variety in a quasiaffine space must be a point).  If we had a
fine moduli space for the varieties in question this would be
automatic, and in this case the family would be trivial and not just
isotrivial. In general we only have a coarse moduli space.  However,
it is well known that,
except for Enriques surfaces, this can be made fine by introducing a
level $n$ structure ($n\ge 3$) on $H^2(\cdot,{\mathbb Z})$. More
precisely we must trivialize the local system whose fibers are
$H^2(\cdot,{\mathbb Z}/3{\mathbb Z})$, and this can be done by taking
a finite \'etale cover of the base space.  This rigidification of
$H^2(\cdot,{\mathbb Z}/3{\mathbb Z})$ provides a fine moduli space
because any element of finite order of $GL_N({\mathbb Z})$ which is
the identity modulo $n\ge 3$ is the identity.  (For Enriques surfaces
we must instead rigidify the cohomology with coefficients in ${\mathbb
Z}/3{\mathbb Z}$ of its K3 cover.)  So to show that compact families
are isotrivial we just have to show that the corresponding coarse
moduli spaces are quasiaffine.

The moduli space of
hyperelliptic surfaces is affine (\cite{BPV}) so theorem 1.1 follows
immediately for these surfaces because the moduli space cannot contain
complete subvarieties of positive dimension,
and hence all fibers in a family are isomorphic and hence
the family is isotrivial by the remarks above. Similarly by theorem 1 of
\cite{RB2} the moduli space of Enriques surfaces is quasiaffine so
theorem 1.1 again follows.  For abelian surfaces theorem 1.1 follows
from the case of K3 surfaces by looking at the associated Kummer
variety and using the fact that the Picard number of
a Kummer surface is 16 plus the Picard number of the abelian surface.
For families of non-projective K3 surfaces, theorem 1.1
follows from a result of Fujiki \cite[Theorem 4.8 (1)]{F}.  So theorem
1.1 follows from the case of projective K3 surfaces, which we treat at
the end of section 2. The main idea of the proof is to construct
automorphic forms on the period space of marked K3 surfaces with
Picard lattice containing a given lattice, such that all zeros of the
automorphic form correspond to K3 surfaces with a larger Picard
lattice. Then we use the fact that the zeros of an automorphic form
give an ample divisor on the moduli space.

We construct these automorphic forms using the denominator 
function of the fake monster Lie algebra
described in  \cite{RB1}, in much the same way that the denominator formula
of the fake monster Lie superalgebra was used in 
\cite{RB2} to prove that the moduli space of Enriques surfaces
is quasiaffine. We can get the following more precise results by
looking more carefully at the automorphic forms we construct.  Recall
that an $S$-K3 surface $X$ for some Lorentzian lattice $S\subset II_{3,19}$
is a K3 surface with a fixed primitive embedding of $S$ into the Picard
group such that the image of $S$ contains a semi-ample class.
(A  semi-ample class is a class $D$ such that $D^2>0$ and 
$D.C\ge 0$ for all curves $C$ on the K3 surface $X$.)

\begin{theo}
There is an automorphic form on the period space of marked $S$-K3
surfaces which vanishes only on divisors $t^\perp$ of vectors $t$ in
the dual $T'$ of $T=S^\perp$ with $0>(t,t)\ge -2$.  (The period space
is an open subset of the space of 2 dimensional positive definite
subspaces of $T\otimes \mathbb{R}$, and $t^\perp $ is the set of
subspaces in the period space orthogonal to $t$; see section 2.)
\end{theo}

By taking $S$ to be a 1 dimensional lattice generated by a vector of
norm $2$ we get: 

\begin{theo} Any family  $f:X \to B$ of smooth,
 polarized K3 surfaces with polarization of degree 2 over a projective
variety $B$ is isotrivial.

\end{theo}

If the period point of a marked $S$-K3 surface lies on the zero locus
of the automorphic form of theorem 1.2 then the surface is $S$-bad in
the following sense.

\begin{defi}  We will call an  $S$-K3 surface    $S$-bad if $S$ is
contained
in a sublattice $S_1$ of the Picard lattice 
such that $\dim(S_1) = \dim(S)+1$ and $|\det(S_1)| \le 2|\det(S)|$.

\end{defi}

Here are some examples of $S$-bad K3 surfaces. If there is a norm $-2$
vector in the Picard lattice orthogonal to $S$ then the K3 surface is
obviously $S$-bad (take $S_1$ to be the lattice generated by $S$ and
this vector). Any K3 surface whose period point lies on one of the
zeros of the automorphic form of theorem 1.2 is also $S$-bad (take
$S_1$ to be generated by $S$ and a vector of the Picard lattice whose
projection to $T'$ is $t$). If $S$ is generated by a vector of norm
$2n$, so that $S$-K3 surfaces include K3 surfaces with a polarization of
degree $2n$, then it is easy to check that a K3 surface is $S$-bad if
its Picard lattice has a vector $D$ of degree $k$ such that $-2\le
(D,D)-k^2/2n<0$. (By adding a multiple of the polarization vector to
$D$ and possibly changing the sign of $D$ we can also assume that
$0\le k\le n$.)

As an immediate consequence of theorem 1.2  we get: 

\begin{theo} Any non isotrivial family  $f:X \to B$ of  $S$-K3 surfaces
over a projective variety $B$ has at least
one $S$-bad fiber.
\end{theo}

We now show why most of the conditions of theorem 1.1 are necessary.

If we take a fixed K3 surface of Picard number $n$ and blow it up in 
a varying point we get a complete non isotrivial family of
non minimal surfaces of Picard number $n+1$. Hence theorem 1.1
does not extend to non minimal surfaces. 

We show that the K\"ahler condition is necessary.  Take two distinct
elliptic curves $A$ and $B$ and let $B_{A}$ be the sheaf of germs of
analytic maps from $A$ to $B$. The sheaf $B_{A}$ is a sheaf of abelian
groups on $A$ and $H^{1}(A, B_{A})$ classifies all primary Kodaira
surfaces with base $A$ and fiber isomorphic to $B$ (\cite{BPV}, V.5.,
p. 143-147).  Represent $B$ as a quotient ${\mathbb
C}/\Lambda$. Choose a nonzero element $c \in H^{2}(A,\Lambda)$ and let
$T \subset H^{1}(A, B_{A})$ be the preimage of $c$ under the natural
map coming from the short exact sequence $0 \to \Lambda \to {\mathcal
O}_{A} \to B_{A} \to 0$. The coset $T$ is isomorphic to the quotient
(stack) $H^{1}(A,{\mathcal O}_{A})/H^{1}(A,\Lambda) = A/{\mathbb
Z}^{2}$. By pulling back to $A$ the universal family over $T$ we
obtain a complete non-isotrivial family of primary Kodaira
surfaces. The latter are never K\"{a}hler since their first Betti
number is equal to 3 and thus theorem 1.1 does not hold for
non-K\"{a}hler surfaces.

The moduli space of curves
of genus greater than 2 
contains compact curves, so the analogue
of theorem 1.1 for curves of high genus is false. By taking products of 
these families with curves of genus 0,1, or greater
than 1 we can find compact non isotrivial  families of minimal surfaces
of constant Picard number 2 with Kodaira dimensions $-\infty$, 1, or 2. 
So theorem 1.1 is not true for surfaces of nonzero Kodaira dimension.

Over fields of characteristic $p>0$ there are complete
non isotrivial families of supersingular K3 surfaces
of constant Picard number 22. (See \cite{S}.) So theorem 1.1 is false
in nonzero characteristics (with ``K\"ahler'' replaced by 
``projective''). However  a plausible analogue 
of the characteristic zero result might be that in a complete
nonisotrivial family of supersingular K3 surfaces the Artin invariant
must jump.

In section 3 of the paper we construct a non isotrivial family of
smooth polarized K3 surfaces.  By
theorem 1.1 the family cannot have constant Picard number and we find
some explicit examples of  $S$-bad fibers for some 1-dimensional $S$.
This shows that the condition about constant Picard number in theorem
1.1 cannot be omitted in the case of K3 surfaces. (The hypothesis
about constant Picard number is irrelevant 
for Enriques or hyperelliptic surfaces because all such surfaces
have Picard numbers 10 and 2 respectively.)  

The K3 surfaces of the family in section 3 all have a fixed point free
involution with an Enriques surface as the quotient. We use
theorem 1.1 to show that there is no global fixed point free
involution of the whole family.

The only condition in theorem 1.1 that we have not shown is necessary 
is the restriction to 2 dimensional varieties. We do not know
what happens in higher dimensions.

\noindent
{\bf Acknowledgments:}
We would like to thank  J. Jorgensen 
and A. Todorov for some interesting discussions about the results
of this paper. In particular R. E. Borcherds would like
to thank A. Todorov for explaining their idea of using
automorphic forms to attack the Shafarevich conjecture for K3 surfaces;
this led to the idea  of using the automorphic forms  of
[B95] to prove results about moduli spaces of surfaces.
Jorgenson and Todorov independently noticed that  results about 
isotriviality of families of K3 surfaces would follow from the existence of
automorphic forms with known zeros, and their preprint
\cite{JT96} has some overlap with the results of this paper, 
such as theorem 1.3. (The methods used in \cite{JT96} to construct
automorphic forms are quite different from the methods in this paper
and are based on regularized determinants as in \cite{JT94}.)
We would also like to thank the referee for pointing out some gaps.   

\section{Automorphic forms on moduli spaces of K3 surfaces} 

In this section we construct some automorphic forms with known zeros
on certain period spaces of marked K3 surfaces with extra structure,
which gives explicit examples of ample divisors on the corresponding
moduli spaces since the set of zeros of an automorphic form is an
ample divisor.  The extra structure consists of a fixed primitive
embedding of some lattice $S$ of signature $(1,m)$ in the Picard
lattice of the K3 surface.  We call such a K3 surface an $S$-K3
surface.  We regard $S$ as a fixed sublattice of the lattice
$II_{3,19}$, and write $T$ for the lattice $S^\perp$ of signature
$(2,19-m)$.  The (hermitean) symmetric space of the lattice $T$ is the
set of norm 0 points in the complex projective space of
$T\otimes{\mathbb C}$ whose real and imaginary parts span a 2
dimensional positive definite subspace of $T\otimes {\mathbb R}$.  We
recall from \cite{BOV} that the moduli space of $S$-K3 surfaces can be
identified with the quotient of the symmetric space of the lattice
$T=S^\perp$ by some arithmetic group.  The Baily-Borel theorem implies
that the zero locus of an automorphic form is an ample divisor on a
compactification of the moduli space, and hence the 
set of points of the quotient where an automorphic form does not vanish
is a quasiaffine variety.

The main result of this section is a proof of theorem 1.2, which
states that there is an automorphic form on the space of marked $S$-K3
surfaces which vanishes only on divisors of vectors $t\in T'$ with
$0>(t,t)\ge -2$.
We will construct this automorphic form  by first embedding the lattice
$T$
in the lattice $II_{2,26}$ and then restricting a certain automorphic form 
of weight 12 for this lattice to the symmetric space of $T$.

We can find a geometric interpretation of the K3 surfaces whose
period points lie on the divisors in theorem 1.2 as follows.  We know
that the Picard lattice contains $S$. As $t\in T'$ and $II_{3,19}$ is
unimodular we can find a vector $D\in II_{3,19}$ whose projection into
$T$ is $t$. Then the lattice $\langle S,D\rangle$ generated by $S$ 
and $D$ has the properties

\begin{itemize}
\item  $\langle S,D\rangle$ has signature $(1,m+1)$

\item $|\det(\langle S,D\rangle)|\le 2|\det (S)|$.

\end{itemize}
because the projection of $v$ into the orthogonal complement of $S$
has norm of absolute value at most 2.  Hence the Picard lattice of the
K3 surface contains a lattice with the properties above.
In particular any  K3 surface for which there
is a norm $-2$ vector in $S^\perp$ satisfies the condition above,
as we can take $D$ to be this norm $-2$ vector.

We now prove theorem 1.2.  

{\bf Proof.} We first construct some primitive embeddings of $T$ into
$II_{2,26}$.  Corollary 1.12.3 of Nikulin \cite{NIK2} implies that we
can primitively embed any lattice $T$ into the unimodular lattice
$II_{2,26}$ provided that $T\otimes \mathbb{R}$ embeds into
$II_{2,26}\otimes \mathbb{R}$ and the minimum number of generators of
$T'/T$ is less than $\dim(II_{2,26})-\dim(T)$.  We can therefore find
a primitive embedding of our lattice $T$ into $II_{2,26}$ because the
rank of the group $T'/T$ is at most the dimension of $S$, so this rank
plus the dimension of $T$ is less than the dimension of
$II_{2,26}$. We will write $U$ for the orthogonal complement $T^\perp$
of $T$ in $II_{2,26}$.  Then $T$ and $U$ have the same determinant as
$T$ is a primitive sublattice of $II_{2,26}$.

We recall some properties of the function $\Phi$ defined in example 
2 of section 10 of \cite{RB1}. The properties of $\Phi$ we will use
are that $\Phi$ is an automorphic form on the hermitian symmetric 
space of $II_{2,26}$ and its only zeros lie on the divisors of norm 
$-2$ vectors of $II_{2,26}$. Some other properties of $\Phi$ which 
we will not use are that its zeros all have multiplicity 1, it has weight 12, 
it is the denominator function of the fake monster Lie algebra, 
its Fourier series is explicitly known, and it can be written explicitly
as an infinite product.

The restriction of $\Phi$ to the hermitian symmetric space of $T$ is
an automorphic form, but will be identically 0 whenever $U$ contains a
norm $-2$ vector. We can get around this by first dividing $\Phi$ by a
product of linear functions vanishing on the divisors of each of these
norm $-2$ vectors before restricting it.  This restriction is an
automorphic form as in pages 200-201 of \cite{RB1}. So in all cases we
get an automorphic form $\Phi_T$ on the hermitian symmetric space of
$T$ whose only zeros lie on the hyperplanes of norm $-2$ vectors of
$II_{2,26}$.  Although we do not need it we can work out the weight of
$\Phi_T$ as follows: the weight is increased by 1 each time we divide
$\Phi$ by a linear function, so the final weight is the weight (=12)
of $\Phi$ plus half the number of norm $-2$ vectors of $U$.

We would like to know the zeros of $\Phi_T$ in terms of vectors of $T$
rather than in terms of norm $-2$ vectors $r$ of the larger lattice
$II_{2,26}$. These zeros correspond to the hyperplanes of the negative
norm projections of the norm $-2$ vectors $r$ of $II_{2,26}$ into $T$.
The projection of $r$ into $U$ has norm at most 0 as $U$ is negative
definite, so the projection $t$ of $r$ into $T$ is a vector of $T'$
with $0>(t,t)\ge -2$.  This proves theorem 1.2.

For the period space of marked Enriques surfaces there is an
automorphic form vanishing exactly on the points orthogonal to $-2$
vectors \cite{RB2}, so it is natural to ask if there is an automorphic
form for polarized K3 surfaces vanishing exactly on the points
orthogonal to a norm $-2$ vector in $S^\perp=T$, which would be much
stronger than the result above. Theorem 1.3 says there is such a form
for K3 surfaces with a polarization of degree 2, but Nikulin
\cite{NIK1} has shown that no such form can exist for some large
values of the polarization.

The zeros of the form in theorem 1.2 do not always have multiplicity
one; in fact they often have some zeros of high multiplicity. We can
work out the multiplicity of the zeros by counting numbers of vectors
in the dual $U'$ of the lattice $U$ with given norm and given image in
$U'/U$. (But notice that some hyperplanes can have higher multiplicity
than one might expect because they get zeros from more than one vector
$t$.)

We now give some examples for polarized K3 surfaces, so we
take $S$ to be a one dimensional lattice spanned by 
a primitive vector of norm $2n$ for some positive  integer $n$. 
We can parameterize embeddings of $T$ into $II_{2,26}$
by primitive norm $-2n$ vectors $v$ in $-E_8$.  To do
this we simply identify $T=(-2n)\oplus(-E_8)\oplus (-E_8)\oplus
H\oplus H$ with the sublattice $\mathbb{Z} v\oplus(-E_8)\oplus (-E_8)\oplus
H\oplus H$ of $II_{2,26}=(-E_8)\oplus(-E_8)\oplus (-E_8)\oplus H\oplus
H$. The lattice $U$ is then the orthogonal complement of $v$ in $-E_8$.

Here is a table of the number of norm $-2$ roots of $U$ and the
numbers of vectors $a$ in the lattice $U'$ of norm
greater than $-2$ for values of $k=(a,v)$ between 0 and $n$.

$ 
\begin{array}{cccccccccc}
     2n &roots&k=0&k=1&k=2&k=3&k=4&k=5&k=6&k=7\\
 2&126&1&56\\
 4&84&1&64&14\\
 6&74&1&54&27&2\\
 8&126&1&0&56&0&1\\
 8&56&1&56&28&8&0\\
10&60&1&44&33&12&1&0\\
12&46&1&48&30&16&3&48&10\\
14&44&1&42&35&14&7&0&21&2\\
14&72&1&28&27&27&1&1&27&0\\
     2n &roots&k=0&k=1&k=2&k=3&k=4&k=5&k=6&k=7

\end{array}
$

The numbers of vectors for other values of $k$ can be worked 
out using the fact that this number does not change if $k$ is replaced
by $2n+k$ or by $-k$. For some values of $n$ there is more than 
one line because the lattice $E_8$ can have several orbits of vectors of
the same norm, corresponding to several different automorphic forms.
The first line for $2n=8$ corresponds to
a non primitive norm $8$ vector of $E_8$ so does not correspond
to a primitive  embedding of $T$ into $II_{2,26}$. 
Some of the entries are 0, corresponding to the fact that the automorphic
forms do not always vanish on all the divisors of theorem 1.2
(so the divisor in theorem 1.2 is not necessarily a minimal 
ample divisor). 

{\bf  Example 2.1.} 
We will work out exactly what the automorphic form for $2n=2$
looks like. Its weight is (weight of $\Phi$)+(number of roots of
$E_7)/2=12+126/2=75$. The zeros of this form come by taking $k=0$ or 1
in theorem 1.2. For $k=0$ we get a contribution of 1 to the multiplicity
of the divisor of each norm $-2$ vector in $T$. For $k=1$ we get a
contribution of 56 for each norm $-1/2$ vector in the dual $T'$ of
$T$. This does not mean that the automorphic form has zeros of
multiplicity 56, because twice a norm $-1/2$ vector of $T'$ is a norm
$-2$ vector of $T$ which has even inner product with all vectors of
$T$. In particular the divisor of the norm $-2$ vectors of $T$ is
reducible: it has two components $E_1$ and $E_2$ corresponding to norm
$-2$ vectors which have odd inner product with some vector of $T$ and
to norm $-2$ vectors which have even inner product with all vectors of
$T$. The divisor $E_1$ is a zero of the automorphic form of
multiplicity 1, and the divisor $E_2$ is a zero of multiplicity
$1+56=57$, and these are all the zeros. 

In particular this example proves theorem 1.3 of the introduction.

{\bf Example 2.2.} Nikulin conjectured in \cite{NIK1} that there are
only a finite number of lattices $S$ such that there is an automorphic
form vanishing only on $S$-K3 surfaces which have a norm $-2$ vector
in $S^\perp$. We can find a few examples of lattices $S$ with this
property. Firstly if $S$ is a unimodular Lorentzian lattice in $L$
then $T=T'$ so $S$ has this property.  The unimodular Lorentzian
lattices in $L$ are $II_{1,1}$, $II_{1,9}$, and $II_{1,17}$. Secondly
if $S$ has determinant 2 and dimension 1 mod 8 then it has the
property above as in example 1, so we also get the lattices $(2)$,
$(2)\oplus (-E_8)$, and $(2)\oplus (-E_8)\oplus (-E_8)$.

{\bf  Example 2.3.} 
If $2n$ is 4,6,8, or 10 then we can see from the table above
that we can assume that either $(D,D)=-2, (D,P)=0$ or $(D,D)=0$. 
Hence if the period of the K3 surface is on a zero of the automorphic form
then either the surface is singular (in the sense that the Picard
group contains a $-2$ vector orthogonal to $S$) or its Picard lattice
contains a nonzero element with zero self intersection number.

\begin{lemma}
Any family of K3 surfaces with constant Picard number 
is, after a finite \'etale base change, a family of $S$-K3 surfaces
for some Lorentzian lattice $S$.
\end{lemma}

Proof. By the assumption about constant rank, and since the Picard
group of a K3 surface is always a primitive sublattice of $H^2$, the
Picard groups form a sub-local system of the local system of $H^2$'s.
Since the monodromy action on the Picard group is finite we can,
after an \'etale base change, assume that this subsystem is constant. 
This means that there is a primitive embedding of a constant local system 
with fiber $S$ into the local system of $H^2$'s. By definition 
this is a family of $S$-K3 surfaces.
This proves lemma 2.1.

We now prove theorem 1.1.  By the remarks after theorem 1.1 we can
assume that all the surfaces in the family are K3 surfaces. 
By lemma 2.1 we can assume that we have a family of $S$-K3
surfaces for some lattice $S$.  
As the
K3 surfaces are projective, the lattice $S$ is Lorentzian. By theorem
1.2 and the remarks near the beginning of this section, if the family
is not isotrivial there must be surfaces whose Picard lattice is
at least 1 more than the dimension of $S$. This contradicts the fact
that all surfaces in the family have the same Picard number and proves
theorem 1.1.

\section{Some examples}

In this section we construct an example of a complete non isotrivial
family of smooth polarized K3 surface, to show that the hypothesis
about the Picard number in theorem 1.1 cannot be left out. 

Suppose that $A\to C$ is a complete one parameter family of
principally polarized family of abelian surfaces. Such a family 
exists because the boundary of the Satake compactification
of the moduli space has codimension $2>1$.  Set $f: K=A/\{\pm 1\}\to
C$. Since $f$ is isotrivial in a neighborhood of its critical locus,
we can simultaneously resolve the singularities via a map
$\sigma: \tilde K\to K$ to get $\tilde K\to
C$, a family of smooth K3 surfaces. 

Let $\Theta$ be the relative principal polarization on $A$.  It is
well known that $2\Theta$ is the pullback of a relatively ample
divisor $H$ on $K$. Let $E$ denote the exceptional locus of
$\sigma$. Then it is well known that on each geometric fiber
$E$ is uniquely even. So on the geometric generic fiber 
$\tilde K_{\bar\eta}$ there is a unique divisor class $\bar L$ 
such that $E_{\bar \eta} \sim 2\bar L$.
Recall that the \'etale cohomology group 
$H^1_{et}(\cdot,{\mathbb G}_m) $ is the Picard group
$Pic(\cdot)$,
and $H^2_{et}(\cdot,{\mathbb G}_m)$ is the Brauer group
$Br(\cdot)$.
{}From the Hochschild-Serre spectral sequence 
$$
E^{pq}_2=H^p(Gal({\mathbb C}(\bar \eta)/{\mathbb C}(\eta)),
H^q_{et}(\tilde K_{\bar\eta},{\mathbb G}_m))
\Rightarrow 
H^{p+q}_{et}(\tilde K_\eta,{\mathbb G}_m)
$$
we get an exact sequence 
$$
0\rightarrow Pic(\tilde K_\eta)\rightarrow Pic(\tilde
K_{\bar\eta})^{Gal({\mathbb C}(\bar \eta)/{\mathbb C}(\eta))}
\rightarrow Br({\mathbb C}(\eta)).
$$ Since the Brauer group $Br({\mathbb C}(\eta))$ of the
function field of a complex curve is trivial (Tsen's theorem) it
follows that $E$ is even on the generic fiber $\tilde K_\eta$. It
follows that $E\sim 2L+V$ for some $L$ and for some $V$ supported on
fibers.  By taking a ramified double cover of the base if necessary we
can assume that $V$ is even so that $E$ is even. Say $E\sim 2M$.  Put
$B=2\sigma^*(H)-M$, so that $B$ is a polarization of degree 8 on
$\tilde K$ provided that $H$ is very ample on $K$, which is equivalent
to the abelian surface $A$ being indecomposable as a principally
polarized abelian variety.  (See \cite{GH}, pages 773-787.)  If,
however, $A$ is decomposable and therefore a product of elliptic curves
then $B$ is merely semi ample. 
Nevertheless take $S$ to be the lattice
generated by $B$. We have constructed a complete non isotrivial 
family of $S$-K3 surfaces.
(The divisor $B$ is not ample on every fiber;
if we want this as well we can take the divisor class $D=B+\sigma^*(H)$ 
which provides a polarization of degree
$(3H)^2+M^2 = 3^2\times 4 -8=28$.)

{\bf Example 3.1.}  Suppose that $0\in C$ is such that $A_0$ is
isomorphic to $X\times X'$ for elliptic curves $X$ and $X'$.  Then
there is a morphism $\alpha:K\to {\mathbb P}^1\times{\mathbb P}^1$ of
degree 2, since ${\mathbb P}^1$ is the Kummer variety of an elliptic
curve.  Let $F_1,F_2$ be fibers of the projections ${\mathbb
P}^1\times {\mathbb P}^1\to {\mathbb P}^1$.  Let $D_i$ be the pullback
of $F_i$ to $\tilde K_0$.  We know that $(D_i,D_i)=0$ as $D_i$ is a
fiber of a morphism onto a curve.  Moreover $(D_i,E_0)=0$, where $E_0$
is the exceptional locus on $\tilde K_0$.  Since $B=2(D_1+D_2)-E/2$ it
follows that $(B,D_i)=4$, so that the lattice generated by $S=\langle
B\rangle$ and $D_1$ has discriminant $-16$. This is an example of an
$S$-bad fiber.

{\bf Remark.} Each of the K3 surfaces in the family above is a Kummer
surface and so by a result of \cite{JHK} has a fixed point free
involution (not necessarily unique!) such that the quotient is an
Enriques surface. But by theorem 1.1 we cannot find a global fixed
point free involution acting on the whole family, otherwise we would
get a complete nonisotrivial family of Enriques surfaces.


\begin{thebibliography}{10}

\bibitem[1]{BPV}{\bf W. Barth, C. Peters, A. Van de Ven} 
{``Compact complex surfaces'', Springer Verlag, 1984. }
\bibitem[2]{RB1}{\bf  R. E. Borcherds} 
{\em   Automorphic forms on $O_{s+2,2}(R)$ and infinite products\/},
Inventiones Mathematicae, 120, 1995, pp. 161-213.
\bibitem[3]{RB2}{\bf  R. E. Borcherds} 
{\em   The moduli space of Enriques surfaces and the fake monster 
Lie superalgebra\/}, Topology vol. 35 no. 3, 699-710, 1996.
\bibitem[4]{BOV}{\bf  } {\em   Geometrie des surfaces K3: 
Modules et periodes\/}, Ast\'erisque, 126, 1985.
\bibitem[5]{F}{\bf A.Fujiki,}
{\em On primitively symplectic compact K\"ahler V-manifolds of
dimension four.\/} in 
Classification of algebraic and analytic manifolds, 
(Katata, 1982), 71-250, Progr. Math., 39, 
Birkh\"auser Boston, Boston, Mass., 1983.
\bibitem[6]{GH} {\bf P. Griffiths, J. Harris} ``Principles of
algebraic geometry'', Wiley, 1978. 
\bibitem[7]{JT94}{\bf J. Jorgenson, A. Todorov} An analytic discriminant
for polarized algebraic K3 surfaces, 1994 preprint.
\bibitem[8]{JT96}{\bf J. Jorgenson, A. Todorov} Ample divisors, automorphic
forms, and Shafarevich's conjecture, 1996 preprint. 
\bibitem[9]{JHK}{\bf  Jong Hae Keum} 
{\em   Every algebraic Kummer surface is the $K3$-cover of an 
Enriques surface\/},  Nagoya Math. Journal  118 (1990), pp. 99--110.
\bibitem[10]{NIK1}{\bf  V. Nikulin} 
{\em   The remark  on  discriminants of  K3 surfaces moduli as sets of 
zeros of automorphic forms\/}, Preprint  alg-geom/9512018.
\bibitem[11]{NIK2}{\bf  V. Nikulin} {\em   Integer symmetric bilinear
forms and some of their geometric applications\/}, Izv. Acad. Nauk
SSSR, Ser. Math. 43 (1979), no. 1, 111-177, 238. English translation
in Mathematics of the USSR Izvestia, Vol. 14 No. 1 1980, 103-167.
\bibitem[12]{S}{\bf A. N. Rudakov, T. Zink, I. R. Shafarevich.}
{\em The influence of height on degenerations of algebraic surfaces of
type K3}, Izv. Acad. Nauk SSSR, Ser Math 46, 117-134 (1982),
translation in Math. USSR, Izv. 20 No. 1, 119-135 (1983). Also
reprinted
in  ``Collected mathematical papers'' by I. R. Shafarevich, Springer 1989.

\end{thebibliography}
\end{document}